# The Timing of Bid Placement and Extent of Multiple Bidding: An Empirical Investigation Using eBay Online Auctions

**Sharad Borle, Peter Boatwright and Joseph B. Kadane**

*Abstract.* Online auctions are fast gaining popularity in today's electronic commerce. Relative to offline auctions, there is a greater degree of multiple bidding and late bidding in online auctions, an empirical finding by some recent research. These two behaviors (multiple bidding and late bidding) are of "strategic" importance to online auctions and hence important to investigate.

In this article we empirically measure the distribution of bid timings and the extent of multiple bidding in a large set of online auctions, using bidder experience as a mediating variable. We use data from the popular auction site www.eBay.com to investigate more than 10,000 auctions from 15 consumer product categories. We estimate the distribution of late bidding and multiple bidding, which allows us to place these product categories along a continuum of these metrics (the extent of late bidding and the extent of multiple bidding). Interestingly, the results of the analysis distinguish most of the product categories from one another with respect to these metrics, implying that product categories, after controlling for bidder experience, differ in the extent of multiple bidding and late bidding observed in them.

We also find a nonmonotonic impact of bidder experience on the timing of bid placements. Experienced bidders are "more" active either toward the close of auction or toward the start of auction. The impact of experience on the extent of multiple bidding, though, is monotonic across the auction interval; more experienced bidders tend to indulge "less" in multiple bidding.

Sharad Borle is Assistant Professor of Management, Jones Graduate School of Management, Rice University, Houston, Texas 77005, USA e-mail: sborle@rice.edu. Peter Boatwright is Associate Professor of Marketing, The Tepper School of Business, Carnegie Mellon University, Pittsburgh, Pennsylvania 15213, USA e-mail: pbhb@andrew.cmu.edu. Joseph B. Kadane is Leonard J. Savage University Professor of Statistics and Social Sciences, Carnegie Mellon University, Pittsburgh, Pennsylvania 15213, USA e-mail: kadane@stat.cmu.edu.

## 0. INTRODUCTION

Over the past few years, online auctions have gained immense popularity and perhaps are the most popular form of electronic commerce (Lucking-Reiley, 2000). Most online auction formats are a variant of second-price sealed-bid auction. For example, the auction format followed by eBay is a second-price auction incorporating features of an English outcry auction as well as the second-price sealed-bid auction. The auction has a fixed ending time (a hardclose auction) and it allocates the object to the bidder with the *highest valuation* (the "maximum bid") at a price which is a small increment above the *second highest* "maximum bid." For further details see eBay's auction site, as well as Bajari and Hortaçsu (2003a, b), Ockenfels and Roth (2006) and Roth and Ockenfels (2002).

In this paper, we study the extent of *late* and *multiple bidding* observed in these auctions. "Late" bidding (also referred to as "sniping") in this context is the tendency of a bidder to put in a bid close to the end of an auction, while "multiple bidding" refers to the tendency of bidders to revise their bids as the auction progresses. Late and multiple bidding is a common occurrence and occurs despite the recommendations of the auction site against such behavior (Roth and Ockenfels, 2002).

We investigate distributions of the extent of late bidding and multiple bidding across eBay auctions in a set of 15 consumer product categories after controlling for the effect of bidder experience. We find that most of the categories significantly differ from one another in their extent of late bidding and multiple bidding. We introduce metrics (which we call the "extent of late bidding" and the "extent of multiple bidding") and rank the categories on these metrics. We then argue that the ranking procedure may be a useful step toward building empirical measures of the extent of common/private values in an auction [empirically distinguishing common/private values is an important challenge faced by auction analysts and designers (Armantier, 2002)].

As for the impact of experience (bidder experience) on the extent of late bidding, we find that the impact is *nonmonotonic*, in that higher-experience bidders bid *either earlier or later* in the auction; that is, the start and close of an auction see a greater participation by experienced bidders than the midtime interval of an auction. The impact of experience on multiplicity of bids is in line with that reported earlier in the literature; namely, greater experience reduces the extent of multiple bidding.

## 1. THE DATA

The data were collected from the online auction site www.eBay.com from August 2001 to February 2002 across 27 products in 15 consumer product categories. Bidding histories of all successfully completed auctions (with at least two bids) were recorded. eBay also has "Dutch" auctions in which bidders compete for multiple quantities of a single item and "private" auctions in which the identities and feedback numbers of bidders are not revealed. These auctions were excluded in our analysis. The information for our study was at the level of individual participants and individual auctions, where the specific data we obtained for each auction were the start and end date/time, identifying information for each bidder, each bidder's experience and every bid's date/time. eBay has a point system for every participant which is the net of number of positive and negative feedbacks received by that person. We use this as a proxy for the reputation/experience level of the buyer after appropriately rescaling it to lie in the set of positive integers (Wilcox, 2000). The data were collected by manually querying eBay's completed auctions database and copying data from those pages. A more efficient way of data collection would have been the use of automated agents (popularly known as "web spiders") which greatly increase the efficiency of data collection [in terms of amounts of data collected as well as errors in the collection process (Bapna et al., 2005)]. With the large amounts of freely available web-based data, issues such as collection of "representative data" also become important (see Shmueli, Jank and Bapna, 2005 for a discussion of various practical and methodological issues surrounding sampling web-based data).

Table 1 lists the products across which auction data were collected. It also lists the number of auctions for each product category and the average number of bids and bidders in these auctions. There were a total of 10,144 auctions across these 15 product categories with 55,852 individual bidders participating in these auctions. The category premium wristwatches realized the highest number of bids per auction and the category hair dryer realized the lowest number of bids per auction. These categories respectively also had the largest and the smallest number of bidders per auction. The 15 product categories were chosen so as to have a wide spread of auctions based on our prior subjective judgments on the certainty of private values in these categories (later in the text, in the section on estimated coefficients, we provide a brief discussion of private/common value auctions). The majority of the auctions were of seven-day duration (56% of all auctions) followed by three- and five-day durations (19% and 16%, resp.). The longest duration auctions were for 10 days (about 9% of all auctions), and the minimum duration was two days (four auctions in all).

TABLE 1
*The products for which auction data were collected*

| Product category | Brands | Number of auctions | Mean number of bids per auction | Mean number of bidders per auction |
|---|---|---|---|---|
| Collectible pottery | *Rookwood and Roseville vase* | 1459 | 10.1 | 5.8 |
| Sunglasses | *Gucci, Oakley* | 1372 | 12.1 | 6.3 |
| Golf balls | *Callaway, Titleist* | 1282 | 8.9 | 5.5 |
| Premium wristwatches | *Cartier, Rolex* | 906 | 15.4 | 7.6 |
| Premium writing pens | *Cross, Waterman* | 659 | 6.9 | 4.2 |
| Computer accessories | *Dell 17" non LCD monitor, HP inkjet color printer* | 652 | 8.9 | 5.4 |
| Golf club bags | *Callaway, Ping* | 503 | 12.4 | 6.2 |
| Neckties | *Brioni, Zegna* | 499 | 6.2 | 3.8 |
| Desktop accessories | *Stapler and tape dispenser (any brand)* | 492 | 5.5 | 3.4 |
| Handheld calculators | *Casio, Sharp* | 474 | 6.3 | 3.5 |
| Luggage bags | *American Tourister, Samsonite* | 446 | 7.7 | 4.1 |
| Men's electric shavers | *Any brand* | 421 | 9.6 | 5.4 |
| Electric drills | *Dewalt cordless drill* | 393 | 13.7 | 7.5 |
| Telescopes and microscopes | *Bausch and Lomb microscope, Celestron telescope* | 333 | 12.7 | 6.6 |
| Hair dryer | *Any brand* | 253 | 5.3 | 3.2 |
| MEAN | | **676** | **9.4** | **5.2** |

One of the variables used in our analysis is the *concentration ratio* of bids [(time left)/(total time)], where *total time* is the duration of auction and *time left* is the gap between bid-time and the end-of-auction-time. Figure 1(a) is a histogram plot of this variable across all auctions and all final bids. Further, Figure 1(b) is a histogram plot of the actual *time left* (in hours) for the final bids of bidders in these auctions. For sake of visual clarity Figure 1(b) consists only of seven-day duration auctions; the pattern of the plot, however, does not change if all auctions are included. Almost 17% of all bids are placed in the last hour. About 38% of these (all bids received in the last hour) are received in the final minute, implying that almost 6% of all bids are received in the last one minute of the auction. It is also interesting that the modes of the plot nearly correspond to intervals of 24 hours, which perhaps indicates that particular hours within a day see increased auction activity in terms of bid placement/start-of-auction. This result echoes similar usage regularities observed in use of online search engines (Telang, Boatwright and Mukhopadhyay, 2004). Figure 2 throws further light on the timing of bid placements/start-of-auction. This figure is a histogram plot of the hour of the day the bids are placed and the hour of the day the seller starts his/her auction. As seen from the figure, the most popular time to start an auction is 1900 hr PST (7:00 pm) (the darker columns in the figure). (Note that PST is Pacific Standard Time, USA which is three hours behind Eastern Standard Time, USA.) Further, the time period 6:00 pm to 9:00 pm is the peak time of the day for sellers to start their auctions. This is also the peak time for the bidders to put in their bids (the lighter columns in Figure 2). We do not have information on the geographical location of bidders; it is reasonable to assume that a majority of them are within the United States across all time zones. All sellers, however, are based in the United States.

Figure 3 is a histogram plot of the number of bids per bidder in an auction. This is an indication of the extent of multiple bidding observed across these auctions; an average bidder across all categories places approximately two bids per auction. A majority of the bidders (63% of all bidders) place just a single bid per auction; however, 10% of them place four or more bids per auction.

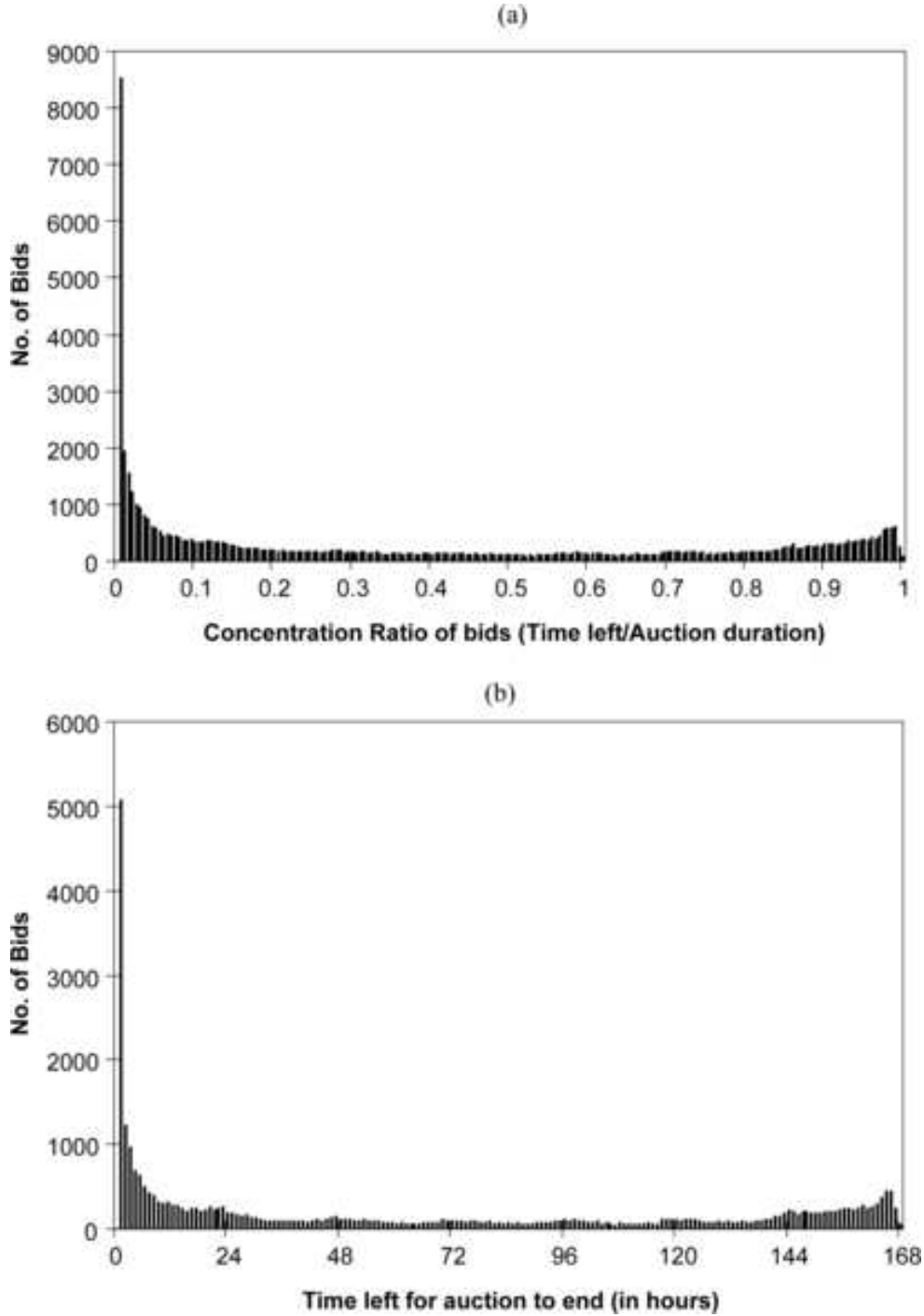

FIG. 1. (a) *Concentration ratio of bids (time left/auction duration) (for the final bids of bidders in all auctions).* (b) *Time gap between bid placement and end of auction (for the final bids of bidders in all seven day auctions).*

Our measure of bidder experience is based on eBay's feedback point system which is the net of number of positive and negative feedbacks received to date on each participant. We use this feedback measure as a proxy for the reputation/experience level of the buyer after appropriately rescaling it to lie in the

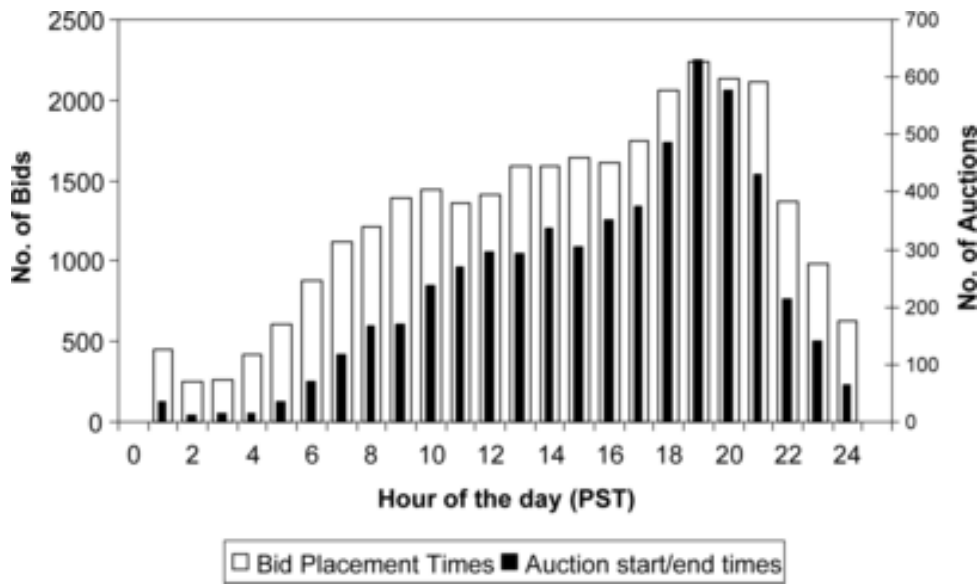

F𝗂𝗀. 2.   *Hour of the day of bid placement and start of auction (for the final bids of bidders in all auctions).*

set of positive integers (Wilcox, 2000). In eBay's system, it is buyers and sellers that receive feedback, so our experience measure reflects the propensity of individuals to complete auctions, either as buyer or seller. An alternative measure of experience, although unavailable to us, would be a total count of the number of bids that each bidder has placed in a specific time period. The second measure reflects participation, regardless of the nature of that participation. Both measures would be correlated with experience and could serve as reasonable proxies. One drawback of the participation measure is that participation does not fully equate with the notion of "experience." In the end, we use the measure that is available to us, which also happens to be that which has been used in the literature (Wilcox, 2000).

The relationship between bidder experience and the timing of bid placement (of the final bid of a bidder) is shown in Figure 4. On the horizontal axis is the concentration ratio of bids [(time left)/(total time)], where *total time* is the duration of auction and *time left* is the time gap between the bid-time and the end-of-auction-time. On the vertical axis is the bidder experience level. For visual clarity experience level is plotted in a logarithmic scale after appropriately scaling it and a random sample of 10,000 bids is plotted instead of all the bids.

An interesting feature of this plot is the nonmonotonic effect of bidder experience on the timing of bid placements (concentration ratio). One can notice increased "clustering" close to the beginning (right-hand side of the plot) and close to the end of auction (left-hand side of the plot), as compared to the middle interval of the auction. Any analysis of such data should allow for such an effect of bidder experience on the timing of bids. In the next section we introduce the distributions of bid timing and the extent of multiple bidding across all 15 product categories.

## 2. THE MODEL

The approach followed in the estimation is to specify different distributions for the "timing of bids"

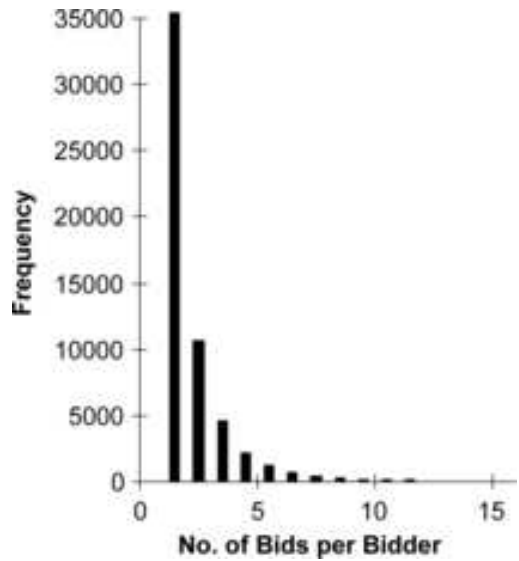

F𝗂𝗀. 3.   *Bids per person.*

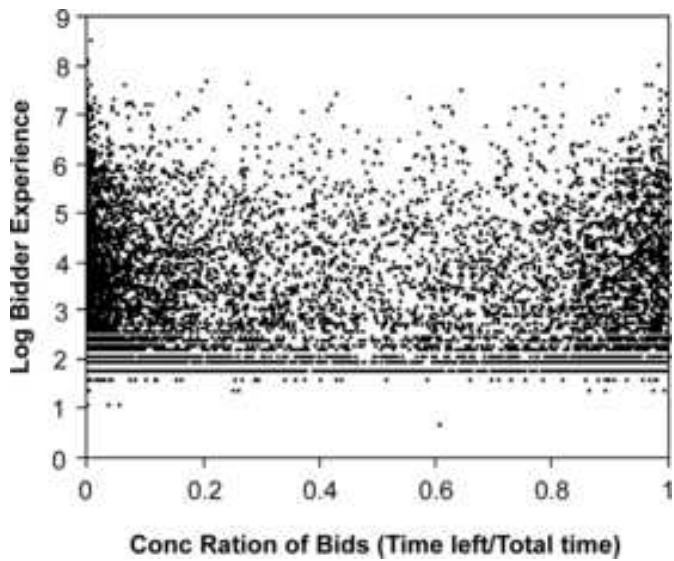

Fig. 4. *Bidder experience and concentration ratio of bids (for the final bids of bidders in all auctions).*

and the "extent of multiple bids" and introduce category-specific parameters across these two distributions. These parameters then are used as metrics to rank-order the product categories based on the extent of late bidding and the extent of multiple bidding observed in them.

### 2.1 Model 1 (the Timing of Bids)

Let $R_{cik}$ be the concentration ratio of bids [(time left)/(total time)], where *total time* is the duration of auction $i$ in product category $c$ and *time left* is the gap between bid-time and end-of-auction-time for the last bid of bidder $k$ in auction $i$ for product category $c$ [we consider only the final bid placed by the bidder; inclusion of all the bids, however, does not substantially change the results].

Thus $R_{cik}$ has the domain $(0,1)$; we assume $R_{cik}$ to have a beta distribution as follows:

$$(1) \qquad R_{cik} \sim Beta(\alpha_{cik}, \beta_{cik}),$$

$$(2a) \qquad \log \alpha_{cik} = \eta + \theta \log EXPB_{cik},$$

$$(2b) \qquad \log \beta_{cik} = \psi'_c + \delta \log EXPB_{cik},$$

where $\alpha_{cik}$ and $\beta_{cik}$ are parameters for the beta distribution. Increase in the value $\alpha_{cik}(\beta_{cik})$ relative to $\beta_{cik}(\alpha_{cik})$ causes the bulk of the distribution to move toward the right (left) side [lower (higher) incidence of late bidding]. When both these parameters decrease (increase), the bulk of the distribution moves away from (toward) the two ends of the distribution. The variable $EXPB_{cik}$ is the experience level of a bidder $k$, in auction $i$, in product category $c$. The parameter $\psi'_c$ is the category-specific metric used to rank-order the product categories based on the extent of late bidding observed in them. Higher value of this parameter for a category would indicate greater extent of late bidding in that category as compared to other categories.

Equations (1) and (2a) suggest that the concentration ratio of bids is a function of the baseline parameter as specified by $\psi'_c$ (2b) (interpreted as a relative measure of the extent of late bidding across auctions in various product categories). Further, this baseline propensity is moderated by the effect of experience on the part of the bidders (as specified by $\delta$ and $\theta$). A positive value for $\delta$ ($\theta$) would imply greater incidence of "late" ("early") bidding on the part of "more" experienced bidders, while a negative value would imply the opposite. A positive (negative) value for *both* $\delta$ and $\theta$ would imply that with greater experience bidders tend to put in their bids away from (close to) the beginning or the end of an auction.

### 2.2 Model 2 (the Extent of Multiple Bidding)

Here we use a measure of the extent of multiple bidding ($BIDS_{cik}$) as the dependent variable. This is the total number of bids placed by a bidder $k$ in auction $i$ in product category $c$. It is modeled as a COM-Poisson (Conway–Maxwell–Poisson, CMP) process (Conway and Maxwell, [1961]),

$$(3) \qquad BIDS_{cik} \sim CMP(\lambda_{cik}, \nu),$$

$$(4) \qquad \log \lambda_{cik} = \psi_c'' + \gamma \log EXPB_{cik},$$

where $\lambda_{cik}$ and $\nu$ are parameters of the COM-Poisson distribution. The COM–Poisson is a generalization of the Poisson distribution with an extra parameter $\nu$ that governs decay. It is a distribution with thicker or thinner tails than the Poisson (Boatwright, Borle and Kadane, [2003]; Shmueli et al., [2005]). The Poisson, the geometric and the Bernoulli distributions are special cases of the COM–Poisson (with $\nu = 1, 0$ and $\infty$, resp.).

The parameter $\lambda_{cik}$, given $\nu$, is proportional to the expected number of bids per bidder [from the probability mass function of the COM-Poisson distribution (Shmueli et al., [2005])]; a "higher" value for this parameter implies a "greater" extent of multiple bidding in the auction.

As in model [(1)], the setup in [(3)] and [(4)] suggests that multiplicity of bids is a function of the baseline metric $\psi_c''$ (the extent of multiple bidding) [[(4)]]. Further, this baseline propensity is moderated by the effect of experience on the part of the bidders as specified by $\gamma$. A positive value for this parameter would suggest that more experienced bidders tend to indulge in a greater incidence of multiple bidding compared to less experienced bidders.

### 2.3 Model Estimation

The Bayesian specification across the two models [[(1)]–[(4)]] is completed by assigning appropriate prior distributions to the parameters to be estimated. The prior distributions and MCMC sampling scheme used in the analysis are given in the Appendix. The models are estimated using a MCMC sampling algorithm, details of which can be obtained from the authors on request.

## 3. THE ESTIMATED COEFFICIENTS

Table [2] reports parameter estimates (posterior mean and standard deviation) not specific to a category while Table [3] reports estimates of parameters specific to a category.

The impact of bidder experience on timing of bid placement is given by $\theta$ and $\delta$ (Table [2]).

Interestingly, both these are estimated to be less than zero, indicating that with increased experience the mass of the beta distribution tends to shift to the two ends. One can interpret this as more experienced bidders tending to be more "active" at the beginning of an auction or toward the close of auction and not in the mid-interval. This nonmonotonic impact of bidder experience on bidding times is an interesting empirical insight not reported earlier. Past work investigating the effect of experience on bid placement times in eBay auctions (Ockenfels and Roth, [2002], [2006]; Wilcox, [2000]) only considered bids received toward the close of auction and concluded that more experienced bidders tend to bid "late" rather than "early." In reporting their empirical results they used a logistic regression with a binary dependent variable (*whether the last bid came in within* 10 *minutes or* 1 *minute of closing*), the experience level being used as an independent variable. What we find is that in general bidders are more active toward the close of auction (low values for the concentration ratio). With increased experience, however, there is a greater tendency on the part of the bidders to bid either toward the close *or* the beginning of auction. For example, in the category *premium writing pens* a high-experience bidder (log experience = 6) is more than twice as likely to bid in the first hour of a seven-day auction as opposed to a low-experience bidder (log experience = 2); the probability increases from 1.1% to 2.7%). Similarly, a high-experience bidder is 1.31 times as likely to bid in the final hour of the auction as opposed to a low-experience bidder (the probability increases from 14.5% to 19.0% ). An analysis failing to allow for a nonmonotonic impact would have led us to conclude that with increased bidder experience a bidder tends to indulge *less* in "late" bidding, a result contrary to past published research. The effect of experienced bidders bidding closer to the auction beginning is the stronger effect and in an analysis not allowing for a nonmonotonic impact outweighs the effect of bidding toward the close of the auction (results of this analysis can be obtained from the authors on request).

TABLE 2
*Model parameter estimates*[*]

| | |
|---|---:|
| $\eta$ | −0.95 |
| | (0.013) |
| $\theta$ | −0.06 |
| | (0.004) |
| $\delta$ | −0.08 |
| | (0.005) |
| $\nu$[**] | 0.00 |
| | (0.000) |
| $\gamma$ | −0.18 |
| | (0.004) |

[*]Shaded cells indicate that the 95% posterior interval for the parameter does not contain 0. The number in parentheses is the posterior standard deviation.

[**]In the estimation in the MCMC $\log(\nu)$ appeared to decrease without bound. Our best approximation of the result was $\nu = 0$, which is the geometric distribution. The model was re-estimated setting $\nu$ to zero.

TABLE 3
*Model parameter estimates*[*]

| | Product category | $\psi_c'$ | $\psi_c''$ |
|---|---|---|---|
| 1 | Collectible pottery | −0.30 (0.023) | −0.23 (0.015) |
| 2 | Sunglasses | −0.31 (0.020) | −0.27 (0.011) |
| 3 | Golf balls | −0.26 (0.021) | −0.48 (0.015) |
| 4 | Premium wristwatches | −0.61 (0.021) | −0.19 (0.012) |
| 5 | Premium writing pens | −0.19 (0.029) | −0.31 (0.021) |
| 6 | Computer accessories | 0.26 (0.028) | −0.41 (0.019) |
| 7 | Golf club bags | −0.41 (0.026) | −0.26 (0.015) |
| 8 | Neckties | −0.18 (0.034) | −0.37 (0.025) |
| 9 | Desktop accessories | 0.03 (0.038) | −0.28 (0.026) |
| 10 | Handheld calculators | 0.12 (0.038) | −0.27 (0.023) |
| 11 | Luggage bags | 0.01 (0.035) | −0.22 (0.020) |
| 12 | Men's electric shavers | −0.04 (0.031) | −0.27 (0.020) |
| 13 | Electric drills | −0.27 (0.027) | −0.28 (0.017) |
| 14 | Telescopes and microscopes | −0.39 (0.032) | −0.18 (0.019) |
| 15 | Hair dryer | 0.03 (0.049) | −0.32 (0.036) |

[*]Shaded cells indicate that the 95% posterior interval for the parameter does not contain 0. The number in parentheses is the posterior standard deviation.

The impact of experience on the extent of multiple bidding is given by the parameter $\gamma$ (Table 2). Earlier studies (Ockenfels and Roth, 2002, 2006; Wilcox, 2000) have reported that more experienced bidders tend to put in fewer bids in the same auction. We also find that more experienced bidders tend to indulge less in multiple bidding as opposed to the less experienced bidders.

The parameters $\psi_c'$ and $\psi_c''$ (the relative extent of late bidding and multiple bidding, resp.) are listed in Table 3 and are also pictorially represented in the box plots of Figure 5(a) and 5(b). To ease interpretation the box plots have been arranged from lower to higher median values, which corresponds to categories with increasing extent of late bidding [Figure 5(a)] and multiple bidding [Figure 5(b)].

An interesting feature of the results is that there is enough variability across categories to differentiate them based on these metrics (the extent of late bidding/multiple bidding). As per the results, the category *premium wristwatches* has the least extent of late bidding while *computer accessories* has the largest extent of late bidding observed among this set of 15 consumer product categories. Similarly, the category *golf balls* has the least extent of multiple bidding, while the category *telescopes and microscopes* has the largest extent of multiple bidding observed.

A relevant question to ask here is does the relative extent of late bidding and multiple bidding in a product category convey some meaningful information? First, the result that after controlling for bidder experience, categories differ in their extent of late bidding and multiple bidding is in itself an interesting empirical result. Second, empirically distinguishing between a private value (PV) and common value (CV) environment is one of the important challenges faced by auction analysts and designers. In a private value auction each bidder values the object differently and knows this valuation before placing the bid. In particular, a bidder's ex post utility from winning the object is not affected by his/her knowledge of other bidders' valuations, since there is no additional information gained by observing the bids of others. In contrast, in a common value auction the object has the same actual value to each bidder ex post, even though different bidders may have access to different information about what that actual value is. Thus, ex ante this value must be estimated by each bidder and these estimates form a distribution around the true value of the object. Past research in this area has primarily used the existence of the winner's curse to distinguish a PV paradigm from a CV paradigm (Paarsch, 1992). In a common value auction the winner is the bidder who has the highest estimate of the "common value" of the object being auctioned and post facto may end up paying more for the object than its true worth; this in short is what is termed the winner's curse (see Milgrom and Weber, 1982 for further explanation of the winner's curse). Armantier (2002) uses the winner's curse to characterize an auction in an experimental setting; Giliberto and Varaiya (1989) use it to characterize FDIC (Federal Deposit Insurance Corporation) auctions of failed banks. Empirically measuring the existence of winner's curse seems to be a well-established method to distinguish the two paradigms. The growth of Internet auctions has spawned further interest in auction theory as applied to these auctions (Lucking-Reiley, 2000). One of the interesting features of these auctions (in comparison to "typical" non-Internet auctions) is the phenomena of late and multiple bidding. As our understanding of the existence of these phenomena and their linkage to the auction environment (a PV or a CV paradigm) grows, we could then use these as additional empirical measures to distinguish the PV and the CV auction paradigms. Based on existing research one might suppose that the extents of late and of multiple bidding may both be indices of common as opposed to private value environment (Hasker, Gonzalez and Sickles, 2004; Bajari and Hortaçsu, 2003b; Roth and Ockenfels, 2002; Unver, 2002; Wilcox, 2000). However, the failure of the rankings in Figure 4 to correlate highly indicates that such an interpretation is questionable (the Spearman rank correlation $-0.496$ is not significantly different from zero).

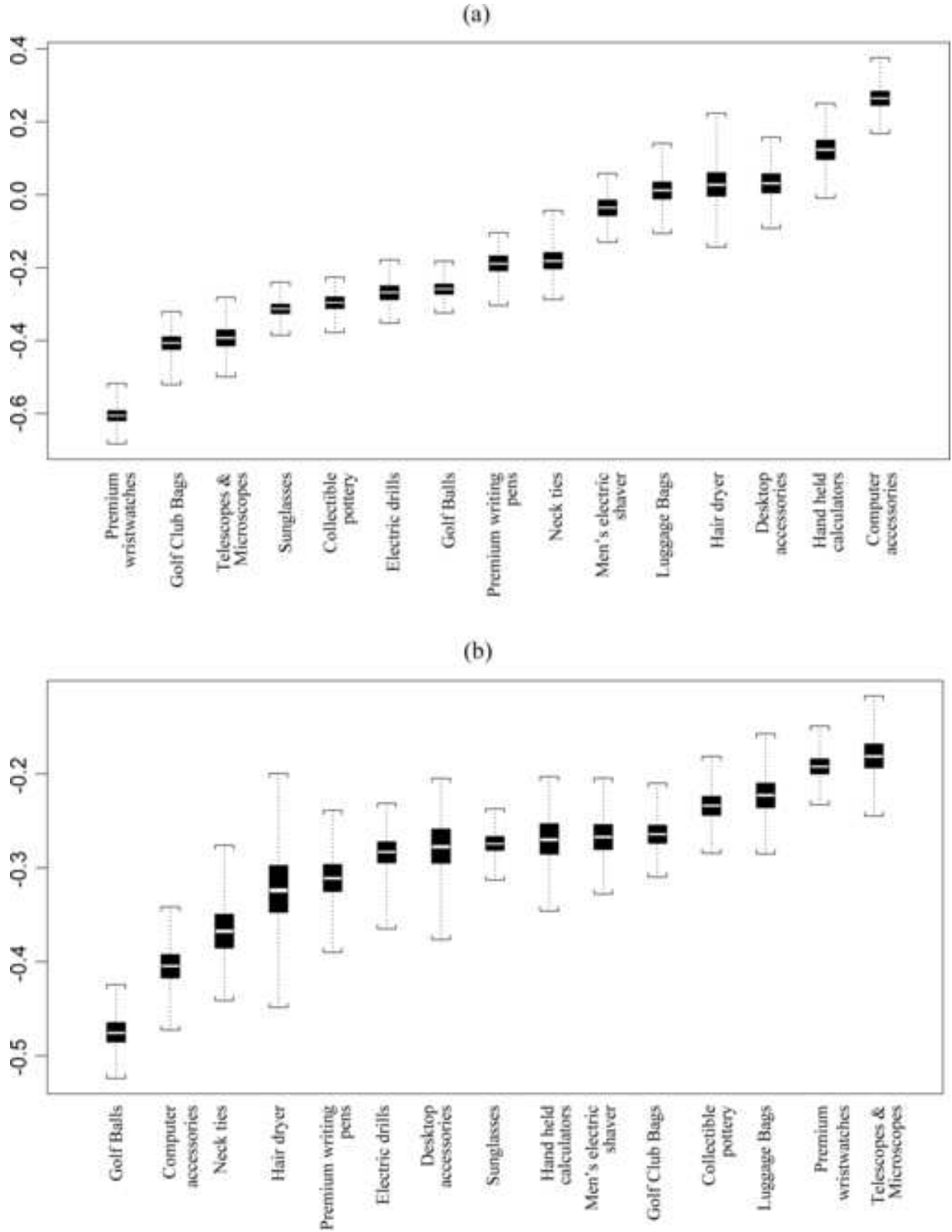

FIG. 5.   *The $\psi_c'$ parameter (extent of late bidding) (higher values indicate "greater" late bidding). The $\psi_c''$ parameter (extent of multiple bidding) (higher values indicate "greater" multiple bidding).*

## 4. DISCUSSION AND CONCLUDING REMARKS

The growing popularity of online auctions has led to much research. The eBay auctions are one of the more prominent among these auctions. Their format is a variant of the second-price sealed-bid auction format. Two additional strategic bidder behaviors (behavior which is not simply due to naive/idiosyncratic bidding but responds to the structure of the auction) observed in these auctions, but not in a second-price sealed-bid auction, are the number of times a bidder bids in the same auction (multiplicity of bids) and the timing (relative to the auction start/end) of his/her bids. These two behaviors and their implications for auction design are of considerable interest to researchers (Bajari and Hortaçsu, 2003b; Hasker, Gonzalez and Sickles, 2004; Ockenfels and Roth, 2002, 2006; Roth and Ockenfels, 2002).

In this paper we investigate over 10,000 auctions across 15 consumer product categories on the eBay auction website www.eBay.com and describe the distribution of the "timing of bids" and the "multiplicity of bids" across these product categories. The Bayesian model used provides a simple and appropriate framework in which to incorporate the structure of the problem and to estimate the relatively large number of parameters. We introduce metrics which summarize the relative extent of late bidding and the relative extent of multiple bidding in an auction. Using these metrics we are able to rank the consumer product categories based on the extent of late bidding and multiple bidding observed in them. Interestingly, we find that most of the product categories vary significantly in the extent of late bidding and multiple bidding observed in them. This variation (as captured by the metrics) can be an important step toward constructing empirical measures of the extent of common/private values in an online auction (Boatwright, Borle and Kadane, 2005). Although most past research implies that the extents of late and of multiple bidding might both be indices of common as opposed to private value bidding in a category, our results question such an interpretation, in that we do not find these measures to be correlated.

We also investigate the effect of bidder experience on the timing of bids and the extent of multiple bids. As reported in earlier studies, we too find that increased experience makes bidders indulge less in multiple bidding (multiple bids placed in the same auction). However, the effect of greater experience (on part of the bidders) on bid placement timings is to move the bids either closer to the auction end time or closer to the auction begin time (a nonmonotonic effect). The more experienced bidders are more active either toward the beginning of the auction or toward the close of the auction as compared to the mid-time interval of the auction. Previous studies on similar data had ignored this nonmonotonic effect and concentrated on the bids received close to the end of auction. This is an interesting empirical finding and a good avenue for further research.

In this work, we restricted our analysis to completed auctions, and our results are thus conditional on auction participation. Incomplete auctions may also serve as an interesting area for future study. There are multiple different auction outcomes that result in incomplete auctions: (a) a seller might withdraw his or her item, (b) there may be no bids at all, and (c) there may be no bids that are high enough to surpass the reservation price. From a statistical modeling viewpoint, incomplete auctions may pose a challenge of sparse data, in that auctions that did not meet the reserve price might have fewer bids overall. From a bidding process point of view, future work could evaluate the extent to which bidders' choices of auctions are influenced by their assessment of the probability that the auction will proceed to completion. In addition, future study could investigate the extent to which such assessments vary by expertise level.

Another potentially interesting study would utilize data where individuals participate in numerous auctions and examine individuals' bidding behavior across auctions as well as possibly the bidders' selection of auctions in which to participate. Future research could also examine how auction bidding is influenced by daily or weekly Internet usage fluctuations. Finally, future work could examine implications of bidder experience for the optimal duration of auctions.

## APPENDIX: MCMC SAMPLING SCHEME AND PRIOR SPECIFICATION

Samplers 1 and 2 are the schematic representation of the MCMC sampling scheme. The priors used in the estimation are as follows: $\delta$ and $\psi'_c$: normal$(0, 100)$; $\theta$ and $\eta$: normal$(0, 100)$; $\nu$: gamma$(2, 1)$; $\gamma$ and $\psi''_c$: normal$(0, 100)$. The first model (the timing of bid placements) corresponding to Sampler 1 is specified by (1), (2a) and (2b). We need to specify our priors over the parameters $\theta, \eta, \delta$ and $\psi'_c$. The parameters $\theta$ and $\delta$ specify the impact of bidder experience on the bid placement times. A priori we do not have specific beliefs about their magnitudes; however, from the empirical distribution (Figure 3)

SAMPLER 1. *The timing of bid placements.*

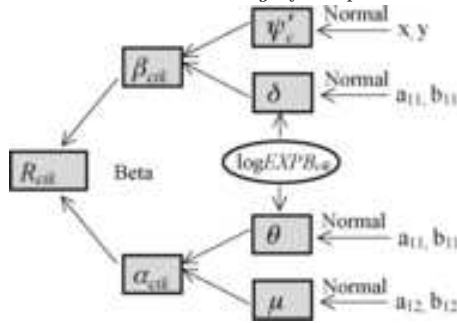

Legend: The variables in boxes are the parameters to be estimated; those in ellipses are the covariates; those without boxes/ellipses are the prior parameters.

SAMPLER 2. *The extent of multiple bidding.*

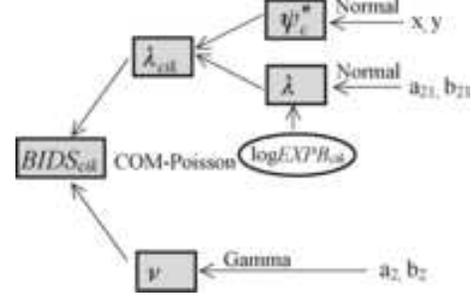

Legend: The variables in boxes are the parameters to be estimated; those in ellipses are the covariates; those without boxes/ellipses are the prior parameters.

one might expect a negative sign on these parameters. A normal$(0, 100)$ is a reasonable representation of this prior belief. The parameter $\psi'_c$ is the relative measure of extent of late bidding in an auction and $\eta$ is the remaining parameter that completes the specification of the beta distribution. Again, a prior of normal$(0, 100)$ represents our lack of specific information about these parameters.

The second model (the extent of multiple bidding) corresponding to Sampler 2 is specified by (3) and (4). The parameter $\nu$ is the decay parameter of the COM-Poisson distribution and is defined over the positive real line. A priori, there is little information on the range of values $\nu$ can take; however, very high values of $\nu$ seem unlikely given the dispersion observed in our data. The prior distribution on $\nu$, a gamma$(2, 1)$ with a mode at 1, is a reasonable representation of our prior belief about the values $\nu$ can take. The parameter $\psi''_c$ is given a prior of normal$(0, 100)$, again representing our lack of specific information about this parameter. The remaining parameter $\gamma$ is the impact of bidder experience on the extent of multiple bidding. As with the corresponding parameters in model 1, a normal$(0, 100)$ represents our prior information about this parameter.

## ACKNOWLEDGMENT

The research of J. B. Kadane was supported by National Science Foundation Grant DMS-98-01401.